\documentclass[11pt]{amsart}
\usepackage{amssymb}
\usepackage{amsmath}
\usepackage{amsthm}
\usepackage{latexsym}
\usepackage{amsfonts}
\usepackage{graphicx}
\usepackage{graphics}
\usepackage[T1]{pbsi}

\newtheorem{lem}{Lemma}
\newtheorem{Def}{Definition}

\theoremstyle{definition}
\newtheorem*{Proof}{Proof}

\newcommand{\dis}{\displaystyle}
\newcommand{\el}{\ell}

\newcommand{\ra}{\;\rightarrow\;}

\newcommand{\bi}{\beta}
\newcommand{\ga}{\gamma }

\newcommand{\OO} {{\varOmega}}

\newcommand{\thi}{\theta }

\newcommand{\la}{\lambda }
\newcommand{\mi}{\mu }

\newcommand{\ti}{\tau }

\newcommand{\oo}{\omega}

\newcommand{\R}{\mathbb{R}}
\newcommand{\Z}{\mathbb{Z}}
\newcommand{\N}{\mathbb{N}}

\newcommand{\intl}{\int\limits}

\newcommand{\cf}{{\mathcal{F}}}

\newcommand{\ct}{{\mathcal{T}}}

\newcommand{\cm}{{\mathcal{M}}}

\newcommand{\ld}{\ldots}

\newcommand{\sm}{\smallsetminus}

 \newcommand{\loc}{\mbox{\footnotesize loc}}
\newcommand{\hs}{\hfill$\square$}
    \newcommand{\vs}{\vspace*{0.2cm} \\}

\begin{document}

\title[DYADIC MAXIMAL OPERATORS]{A sharp integral rearrangement inequality for the dyadic maximal operator and applications}
\author{Eleftherios N. Nikolidakis, Antonios D. Melas}
%
%
\footnotetext{\hspace{-0.5cm}This research has been co-financed by the European Union and Greek national funds
through the Operational Program "Education and Lifelong Learning" of the National Strategic Reference Framework (NSRF).
Aristeia code:MaxBellman 2760, Research code:70/3/11913}
\date{}
\maketitle
\noindent
{\bf Abstract:} We prove a sharp integral inequality for the dyadic maximal operator and give as an application another proof for the computation of its Bellman function of three variables.

Keywords:Bellman, Dyadic, Maximal, Rearrangement.
\noindent
\section{Introduction}\label{sec1}

The dyadic maximal operator on $\R^n$ is defined by
\begin{eqnarray}
\cm_d\phi(x)=\sup\bigg\{\frac{1}{|Q|}\int_Q|\phi(u)|du:\;x\in Q,\;Q\subseteq\R^n\;\text{is a dyadic cube}\bigg\}\hspace*{-1cm}\label{eq1.1}
\end{eqnarray}
for every $\phi\in L^1_{\loc}(\R^n)$, where the dyadic cubes are those formed by the grids $2^{-N}\Z^n$ for $N=0,1,2,\ld\;.$

It is well known that it satisfies the following weak type (1,1) inequality
\begin{eqnarray}
|\{x\in\R^n:\cm_d\phi(x)>\la\}|\le\frac{1}{\la}\int_{\{\cm_d\phi>\la\}}|\phi(u)|du, \label{eq1.2}
\end{eqnarray}
for every $\phi\in L^1(\R^n)$ and every $\la>0$.

Using this inequality it is not difficult to prove the following known as Doob's inequality
\begin{eqnarray}
\|\cm_d\phi\|_p\le\frac{p}{p-1}\|\phi\|_p,  \label{eq1.3}
\end{eqnarray}
for every $p>1$ and $\phi\in L^p(\R^n)$.

It is an immediate result that the weak type inequality (\ref{eq1.2}) is best possible, while (\ref{eq1.3}) is also sharp (see \cite{1}, \cite{2} for general martingales and \cite{16} for dyadic ones).

A way of studying the dyadic maximal operator is by making refinements of the above inequalities.
The above inequalities hold true even in more general settings. More precisely we consider a non-atomic probability space $(X,\mi)$ equipped with a tree structure $\ct$ and define
\[
\cm_\ct\phi(x)=\sup\bigg\{\frac{1}{\mi(I)}\int_I|\phi|d\mi:\;x\in I\in\ct\bigg\}.
\]

Concerning (\ref{eq1.2}) certain refinements have been done in \cite{8} and \cite{9} while for (\ref{eq1.3}) the Bellman function of the dyadic maximal operator has been explicitly computed in \cite{3}.
This is given by
\begin{eqnarray}
B_p(f,F)=\sup\bigg\{\int_X(\cm_\ct\phi)^pd\mi:\;\phi\ge0,\;\int_X\phi d\mi=f,\;\int_X\phi^pd\mi=F\bigg\},  \label{eq1.4}
\end{eqnarray}
for $p>1$ and every $f$ and $F$ such that $0<f^p\le F$.

It is proved in \cite{3} that it equals
\[
B_p(f,F)=F\oo_p(f^p/F)^p, \ \ \text{where} \ \ \oo_p:[0,1]\ra\bigg[1,\frac{p}{p-1}\bigg]
\]
denotes the inverse function $H^{-1}_p$ of $H_p$, which is defined by $H_p(z)=-(p-1)z^p+pz^{p-1}$, for $z\in[1,\frac{p}{p-1}]$.

After this evaluation the second task is to find the exact value of the following function of three variables
\begin{eqnarray}
B_p(f,F,L)=\sup\bigg\{\int_X\max(\cm_\ct\phi,L)^pd\mi:\;\phi\ge0,\;\int_X\phi d\mi=f, \;\int_X\phi^pd\mi=F\bigg\},\hspace*{-1.5cm}  \label{eq1.5}
\end{eqnarray}
for $p>1$, $0<f^p\le F$ and $L\ge f$.

It turns out that
\begin{eqnarray}
B_p(f,F,L)=\left\{\begin{array}{ccc}
                    F\oo_p\Big(\dfrac{pL^{p-1}f-(p-1)L^p}{F}\Big)^p, & \text{if} & L<\dfrac{p}{p-1}f \\ [2.5ex]
                    L^p+\Big(\dfrac{p}{p-1}\Big)^p\Big(F-f^p\Big) & \text{if} & L>\dfrac{p}{p-1}f.
                  \end{array}\right.  \label{eq1.6}
\end{eqnarray}
For this evaluation the author in \cite{3} used the result for (\ref{eq1.4}) on suitable subsets of $X$ and after several calculus arguments he was able to provide a proof of (\ref{eq1.6}).

The Bellman functions have been studied also in \cite{4}. There a more general Bellman function has been computed, namely
\begin{align}
T_{p,G,H}(f,F,k)=\sup\bigg\{&\int_KG(\cm_\ct\phi)d\mi:\;\phi\ge0,\;\int_X\phi d\mi=f,\;\int H(\phi)d\mi=F,\nonumber \\
&K\;\text{measurable subset of}\;X\;\text{with}\;\mi(K)=k\bigg\} \label{eq1.7}
\end{align}
for suitable convex, non-negative, increasing functions $G$ and $H$.
The approach used in \cite{4} is by proving that $T_{p,G,H}(f,F,k)$ equals
\begin{align*}
S_{p,G,H}(f,F,k)=\sup\bigg\{&\int^k_0G\bigg(\frac{1}{t}\int^t_0g\bigg)dt:\;g:(0,1]\ra\R^+\;
\text{non-increasing, continuous} \\
&\text{with}\; \int^1_0g(u)du=f,\;\int^1_0H(g)dt=F\bigg\}.
\end{align*}
The second step then is to evaluate $S_{p,G,H}(f,F,k)$, which in general is a difficult task. Concerning the first step $(T_{p,G,H}=S_{p,G,H})$ the following equality has been proved in \cite{10} stated as\vs
\noindent
{\bf Theorem A.} {\em If $g,h:(0,1]\ra\R^+$ are non-increasing integrable functions and $G:[0,+\infty)\ra[0,+\infty)$ is non-decreasing, then the following is true
\begin{align*}
\sup\bigg\{\int_KG[(\cm_\ct\phi)^\ast]&h(t)dt,\;\phi^\ast=g,\;K\;\text{measurable subset of}\;(0,1]\;\text{with}\;|K|=k\bigg\} \\
&=\int^k_0G\bigg(\frac{1}{t}\int^t_0g(u)du\bigg)h(t)dt.
\end{align*}
}
This can be viewd as a symmetrization principle that immediately yields the equality $T_{p,G,H}=S_{p,G,H}$.

In this paper our aim is to find another proof of (\ref{eq1.6}) by using a variant of Theorem A.
More precisely we will prove the following\vs
\noindent
{\bf Theorem 1.}{\em The following equality is true
\begin{align*}
\sup\bigg\{&\int_KG_1(\cm_\ct\phi)G_2(\phi)d\mi:\phi^\ast=g,\;K\;\text{measurable subset of} \\
&X\;\text{with}\;\mi(K)=k\bigg\}=
\int^k_0G_1\bigg(\frac{1}{t}\int^t_0g\bigg)G_2(g(t))dt,
\end{align*}
where $G_i:[0,+\infty)\ra[0,+\infty)$ are increasing functions for $i=1,2$, while $g:(0,1]\ra\R$ is non-increasing.}

This theorem and some extra effort will enable us to provide a simpler proof of (\ref{eq1.6}).
We also remark that there are several problems in Harmonic Analysis were Bellman functions arise. Such problems (including the dyadic Carleson imbedding theorem and weighted inequalities) are described in \cite{7} (see also \cite{5}, \cite{6}) and also connections to Stochastic Optimal Control are provided, from which it follows that the corresponding Bellman functions satisfy certain nonlinear second-order PDEs. The exact evaluation of a Bellman function is a difficult task which is connected with the deeper structure of the corresponding Harmonic Analysis problem. Until now several Bellman functions have been computed (see \cite{1}, \cite{2}, \cite{3}, \cite{5}, \cite{12}, \cite{13}, \cite{14}, \cite{15}). The exact evaluation of (\ref{eq1.4}) has been also given in \cite{11} by L. Slavin, A. Stokolos and V. Vasyunin which linked the computation of it to solving certain PDEs of the Monge-Amp\`{e}re type and in this way they obtained an alternative proof of the results in \cite{3} for the Bellman functions related to the dyadic maximal operator.

The paper is organized as follows.
In Section \ref{sec2} we give some preliminaries needed for use in the subsequent sections.
In Section \ref{sec3} we prove Theorem 1 while in Section \ref{sec4} we give a proof that the right side of (\ref{eq1.6}) is an upper bound of the quantity: $\int\limits_X\max(\cm_\ct\phi,L)^pd\mi$.
At last in Section \ref{sec5} we prove the sharpness of the above mentioned result.
\section{Preliminaries}\label{sec2}
\noindent

Let $(X,\mi)$ be a non-atomic probability measure space.
\begin{Def}\label{Def2.1}
\ A set $\ct$ of measurable subsets of $X$ will be called a tree if it satisfies the following conditions
\hspace*{-1cm}\begin{enumerate}
\item[i)] $X\in\ct$ and for every $I\in\ct$ we have that $\mi(I)>0$.
\item[ii)] For every $I\in\ct$ there corresponds a finite or countable subset $C(I)\subseteq\ct$ containing at least two elements such that
\begin{itemize}
\item[(a)] the elements of $C(I)$ are pairwise disjoint subsets of $I$
\item[(b)] $I=\cup C(I)$.
\end{itemize}
\item[iii)] $\ct=\bigcup\limits_{m\ge0}\ct_{(m)}$ where $\ct_{(0)}=\{X\}$ and $\ct_{(m+1)}=\bigcup\limits_{I\in\ct_{(m)}}C(I)$.
\item[iv)] We have that $\dis\lim_{m\ra\infty}\dis\sup_{I\in\ct_{(m)}}\mi(I)=0$.
\end{enumerate}
\end{Def}
Examples of trees are given in \cite{3}.\ The most known is the one given by the family of all dyadic subcubes of $[0,1]^n$.
The following has been proved in \cite{3}.
\begin{lem}\label{lem2.1}
For every $I\in\ct$ and every $a$ such that $0<a<1$ there exists a subfamily $\cf(I)\subseteq\ct$ consisting of disjoint subsets of $I$ such that
\[
\mi\bigg(\bigcup_{J\in\cf(I)}J\bigg)=\sum_{J\in\cf(I)}\mi(J)=(1-a)\mi(I).  \]
\end{lem}

We will also need the following fact obtained in \cite{10}.
\begin{lem}\label{lem2.2}
Let $\phi:(X,\mi)\ra\R^+$ and $(A_j)_j$ a measurable partition of $X$ such that $\mi(A_j)>0$ $\forall\;j$.\ Then if $\int\limits_X\phi d\mi=f$ there exists a rearrangement of $\phi$, say $h$  $(h^\ast=\phi^\ast)$ such that $\dfrac{1}{\mi(A_j)}\int\limits_{A_j}hd\mi=f$, for every $j$.
\end{lem}
Here by $\phi^\ast$ we mean the decreasing rearrangement of $\phi$ defined by
\[
\phi^\ast(t)=sup_{e\subset X, |e|=t}[inf_{x\in e}|\phi(x)|], t\in (0,1].
\]

Now given a tree on $(X,\mi)$ we define the associated dyadic maximal operator as follows
\[
\cm_\ct\phi(x)=\sup\bigg\{\frac{1}{\mi(I)}\int_I|\phi|d\mi:\;x\in I\in\ct\bigg\},
\]
for every $\phi\in L^1(X,\mi)$.

We will also need the following well known (see\cite{17}).
\begin{lem}\label{lem2.3}
Let $\phi_1,\phi_2:X\ra\R^+$ be $\mi$-measurable functions. Then the following inequality is always true:
\[
\int_X\phi_1(x)\phi_2(x)d\mi(x)\le\int^1_0\phi^\ast_1(t)\cdot\phi_2^\ast(t)dt
\]
where $\phi^\ast_i$ is decreasing rearrangement of $\phi_i$.
\end{lem}
\section{The rearrangement inequality}\label{sec3}
\noindent

We prove first the following
\setcounter{lem}{0}
\begin{lem}\label{lem3.1}
With the notation of Theorem 1 the following inequality holds
\[
\int_KG_1(\cm_\ct\phi)G_2(\phi)d\mi\le\int^k_0G_1\bigg(\frac{1}{t}\int_0^tg\bigg)G_2(g(t))dt \]
\end{lem}
\begin{Proof}
Following \cite{10} we set
\[
I=\int_KG_1(\cm_\ct\phi)G_2(\phi)d\mi.
\]
Then by using Lemma \ref{lem2.3} we have that:
\[
I\le\int^k_0[G_1(\cm_\ct\phi)/K]^\ast\cdot[G_2(\phi)/K]^\ast dt.
\]
Since $K\subseteq X$ we have that
\[
\begin{array}{ll}
  [G_1(\cm_\ct\phi)/K]^\ast(t)\le[G_1(\cm_\ct\phi)]^\ast(t) & \text{and} \\[1ex]
  [G_2(\phi)/K]^\ast(t)\le[G_2(\phi)]^\ast(t), & \text{for any} \ \ t\in(0,k].
\end{array}
\]
On the other hand, $G_1$ and $G_2$ are increasing functions, therefore
\[
\begin{array}{ll}
  [G_1(\cm_\ct\phi)]^\ast=G_1[(\cm_\ct\phi)^\ast] & \text{and} \\ [1ex]
  [G_2(\phi)]^\ast=G_2(\phi^\ast), &
\end{array}
\]
almost everywhere with respect to the Lesbesgue measure on $(0,k]$. Thus
\[
I\le\int^k_0G_1[(\cm_\ct\phi)^\ast(t)]\cdot G_2(g(t))dt=II.
\]
The last integral now equals
\[
II=\int^k_0G_1[(\cm_\ct\phi)^\ast(t)]dv_2(t),
\]
where $v_2$ is the Borel measure defined on $(0,k]$ by
\[
v_2(A)=\int_AG_2(g(u))du.
\]
Then we have that
\[
II=\int^{+\infty}_{\la=0}v_2(\{t\in(0,k]:(\cm_\ct\phi)^\ast(t)\ge\la\})d \, G_1(\la)=III+IV, \ \ \text{where}
\]
\[
III=\int^f_0v_2((0,k])d\,G_1(\la)=v_2((0,k])[G_1(f)-G_1(0)] \ \ \text{and}
\]
\setcounter{equation}{0}
\begin{eqnarray}
IV=\int^{+\infty}_{\la=f}v_2(\{t\in(0,k]:(\cm_\ct\phi)^\ast(t)\ge\la\})dG_1(\la).  \label{eq3.1}
\end{eqnarray}

Now we will prove that if we set
\[
A_\la=\{t\in(0,k]:(\cm_\ct\phi)^\ast(t)\ge\la\} \ \ \text{and}
\]
\[
\OO_\la=\bigg\{t\in(0,k]:\frac{1}{t}\int^t_0g\ge\la\bigg\},
\]
then $A_\la\subseteq\OO_\la$, for any $\la>f$. Fix such a $\la$.

Since $A_\la$ and $\OO_\la$ are defined in terms of non-increasing functions on $(0,k]$ we must have that
\[
A_\la=(0,|A_\la|], \ \ \text{and} \ \ \OO_\la=(0,|\OO_\la|],
\]
that is they must be intervals with 0 being their common left end-point.
Thus in order to prove that $A_\la\subseteq\OO_\la$ we just need to show that $|A_\la|\le|\OO_\la|$.

For our fixed $\la$ we have that there exists $\bi(\la)\in(0,1]$ such that $\dfrac{1}{\bi(\la)}\dis\int^{\bi(\la)}_0g(u)du=\la$. It's  existence is guaranteed by the fact that $\la>f=\dis\int^1_0g(u)du$. In fact, we can suppose without loss of generality that $g(0^+)=+\infty$, otherwise we work on $\la\in(f,\|g\|_\infty]$. Notice that if $\|g\|_\infty=A$, then $\cm_\ct\phi\le A$ $\mi$-a.e. on $X$.

By the definition of $\OO_\la$ and $\bi(\la)$ it follows that $\OO_\la=(0,\min(\bi(\la),k)]$. Also note that $|A_\la|\le k$. Therefore it suffices to prove that $|A_\la|\le\bi(\la)$. But
\[
A_\la\subseteq\{t\in(0,1]:(\cm_\ct\phi)^\ast(t)\ge\la\}\Rightarrow|A_\la|\le|\{t\in(0,1]:
(\cm_\ct\phi)^\ast(t)\ge\la\}|=\mi(E_\la),
\]
where $E_\la$ is defined by
\[
E_\la=\{x\in X:(\cm_\ct\phi)(x)\ge\la\}.
\]
There exists a pairwise disjoint family of elements of $\ct$, $(I_j)_j$, such that
\begin{eqnarray}
\frac{1}{\mi(I_j)}\int_{I_j}\phi d\mi\ge\la \ \ \text{and} \ \ E_\la=\cup I_j.  \label{eq3.2}
\end{eqnarray}
In fact we just need to consider the family $(I_j)_j$ of elements of $\ct$, maximal under the above integral condition.

By (\ref{eq3.2}) we have that $\dis\int_{I_j}\phi d\mi\ge\la\mi(I_j)$, for any $j$, and so summing the above inequalities with respect to $j$, we conclude that
\[
\int_{E_\la}\phi d\mi\ge\la\mi(E_\la) \ \ \text{or that} \ \ \frac{1}{\mi(E_\la)}\int_{E_\la}\phi d\mi\ge\la.
\]
On the other hand $\bi(\la)$ is defined by the equation:
\[
\frac{1}{\bi(\la)}\int^{\bi(\la)}_0g(u)du=\la.
\]
So we have have the following inequalities
\[
\frac{1}{\mi(E_\la)}\int^{\mi(E_\la)}_0g(u)du\ge\frac{1}{\mi(E_\la)}\int_{E_\la}\phi d\mi\ge\la=\frac{1}{\bi(\la)}\int^{\bi(\la)}_0g(u)du.
\]
implying that $\mi(E_\la)\le\bi(\la)$, since $g$ is non-increasing. Then because of the inequality $|A_\la|\le\mi(E_\la)$ we have $|A_\la|\le|\OO_\la|$. By the above we find that
\[
A_\la\subseteq\OO_\la\Rightarrow v_2(A_\la)\le v_2(\OO_\la).
\]
Now using (\ref{eq3.1}) we get
\[
IV\le\int^{+\infty}_{\la=f}v_2\bigg(\bigg\{t\in(0,k]:\frac{1}{t}\int^t_0g\ge\la\bigg\}\bigg)dG_1(\la), \ \ \text{thus}
\]
\begin{align*}
I&\le\int^{+\infty}_{\la=0}v_2\bigg(\bigg\{t\in(0,k]:\frac{1}{t}\int^t_0g\ge\la\bigg\}\bigg)dG_1(\la)\\
&=\int^k_0G_1\bigg(\frac{1}{t}\int^t_0g\bigg)dv_2(t)=\int^k_0G_1\bigg(\frac{1}{t}\int^t_0g\bigg)G_2(g(t))dt
\end{align*}
by the definition of $v_2$. This completes the proof of Lemma \ref{lem3.1}. \hs
\end{Proof}
We now proceed to the

{\bf Proof of Theorem 1}:
First suppose that $k=1$.
Let $g:(0,1]\ra\R^+$ be a non-increasing function.
We are going to construct a family $(\phi_a)_{a\in(0,1)}$ of functions defined on $(X,\mi)$, each having $g$ as it's decreasing rearrangement $(\phi^\ast_a=g)$, such that
\[
\underset{a\ra 0^+}{\lim\sup}\int_XG_1(\cm_\ct\phi_a)G_2(\phi_a)d\mi\ge\int^1_0
G_1\bigg(\frac{1}{t}\int^t_0g\bigg)G_2(g(t))dt.
\]
Following \cite{10} we let $a\in(0,1)$. Using Lemma \ref{lem2.1} we choose for every $I\in\ct$ a family $\cf(I)\subseteq\ct$ of disjoint subsets of $I$ such that
\begin{eqnarray}
\sum_{J\in\cf(I)}\mi(J)=(1-a)\mi(I).  \label{eq3.2}
\end{eqnarray}
Define $S=S_a$ by induction to be the smallest subset of $\ct$ for which $X\in S$ and for every $I\in S$, $\cf(I)\subseteq S$. We write for $I\in S$, $A_I=I\sm\bigcup\limits_{J\in\cf(I)}J$. Then if $a_I=\mi(A_I)$ we have because of (\ref{eq3.2}) that $a_I=a\mi(I)$.\ It is also clear that
\[
S_a=\bigcup_{m\ge0}S_{a,(m)}, \ \ \text{where} \ \ S_{a,(0)}=\{X\} \ \ \text{and} \ \ S_{a,(m+1)}=\bigcup_{I\in S_{a,(m)}}\cf(I).
\]
We define also for $I\in S_a$, rank$(I)=r(I)$ to be the unique integer $m$ such that $I\in S_{a,(m)}$.
Additionally, we define for every $I\in S_a$ with $r(I)=m$
\[
\ga(I)=\ga_m=\frac{1}{a(1-a)^m}\int^{(1-a)^m}_{(1-a)^{m+1}}g(u)du.
\]
and
\[
b_m(I)=\sum_{S\ni J\subseteq I\atop r(J)=r(I)+m}\mi(J).
\]
We easily then see inductively that
\[
b_m(I)=(1-a)^m\mi(I).
\]
It is also clear that for every $I\in S_a$
\[
I=\bigcup_{S_a\ni J\subseteq I}A_J.
\]
At last we define for every $m$ the measurable subset of $X$, $S_m=\bigcup\limits_{I\in S_{a,(m)}}I$.
Now for each $m\ge0$ we choose $\ti^{(m)}_a:S_m\setminus S_{m+1}\ra\R$ such that
\[
\big[\ti^{(m)}_a\big]^\ast=\Big(g/\big((1-a)^{m+1},(1-a)^m\big]\Big)^\ast.
\]
This is possible since $\mi(S_m\setminus S_{m+1})=\mi(S_m)-\mi(S_{m+1})=b_m(X)-b_{m+1}(X)=(1-a)^m-(1-a)^{m+1}=a(1-a)^m$. It is obvious that $S_m\setminus S_{m+1}=\bigcup\limits_{I\in S_{a,(m)}}A_I$ and that
\[
\int_{S_m\setminus S_{m+1}}\ti^{(m)}_ad\mi=\int^{(1-a)^m}_{(1-a)^{m+1}}g(u)du\Rightarrow\frac{1}
{\mi(S_m\setminus S_{m+1})}\int_{S_m\setminus S_{m+1}}\ti^{(m)}_ad\mi=\ga_m.
\]
Define $\ti_a:X\ra\R^{+}$ by $\ti_a/(S_m\setminus S_{m+1}):=\ti^{(m)}_a$, $m\ge0$.
Using Lemma \ref{lem2.2} we see that there exists a rearrangement of $\ti^{(m)}_a$, called $\phi^{(m)}_a$, for which $\dfrac{1}{a_I}\int\limits_{A_I}\phi^{(m)}_a=\ga_m$, for every $I\in S_{a,(m)}$.
We define $\phi_a:X\ra\R^+$ by $\phi_a(x)=\phi^{(m)}_a(x)$, for $x\in S_m\setminus S_{m+1}$. Clearly $\phi^\ast_a=g$.

Let now $I\in S_{a,(m)}$. Then
\begin{align}
&\frac{1}{\mi(I)}\int_I\phi_ad\mi \nonumber\\
&=\frac{1}{\mi(I)}\sum_{S_a\ni J\subseteq I}\int_{A_J}\phi_ad\mi \nonumber\\
&=\frac{1}{\mi(I)}\sum_{\el\ge0}\sum_{S_a\ni J\subseteq I\atop r(J)=r(I)+\el}\int_{A_J}\phi_ad\mi \nonumber\\
&=\frac{1}{\mi(I)}\sum_{\el\ge0}\sum_{S_a\ni J\subseteq I}\ga_{m+\el}a_J \nonumber\\
&=\frac{1}{\mi(I)}\sum_{\el\ge0}\sum_{S_a\ni J\subseteq I}a\mi(J)\frac{1}{a(1-a)^{m+\el}}\int^{(1-a)^{m+\el}}_{(1-a)^{m+\el+1}}g(u)du \nonumber\\
&=\frac{1}{\mi(I)}\sum_{\el\ge0}\frac{1}{(1-a)^{m+\el}}\int^{(1-a)_{m+\el}}_{(1-a)^{m+\el+1}} \nonumber
g(u)du\cdot\sum_{S_a\ni J\subseteq I\atop r(J)=m+\el}\mi(J) \nonumber\\
&=\frac{1}{\mi(I)}\sum_{\el\ge0}\frac{1}{(1-a)^{m+\el}}\int^{(1-a)^{m+\el}}_{(1-a)^{m+\el+1}}g(u)du\cdot b_\el(I) \nonumber\\
&=\frac{1}{(1-a)^m}\sum_{\el\ge0}\int^{(1-a)^{m+\el}}_{(1-a)^{m+\el+1}}g(u)du \nonumber\\
&=\frac{1}{(1-a)^m}\int^{(1-a)^m}_0g(u)du.  \label{eq3.4}
\end{align}
Now for $x\in S_m\setminus S_{m+1}$, there exists $I\in S_{a,(m)}$ such that $x\in I$ so
\begin{eqnarray}
\cm_\ct(\phi_a)(x)\ge \frac{1}{\mi(I)}\int_I\phi_ad\mi=\frac{1}{(1-a)^m}\int^{(1-a)^m}_0g(u)du=:\thi_m, \label{eq3.5}
\end{eqnarray}
Then for each $a\in(0,1)$ we have that
\begin{align}
\int_XG_1(\cm_\ct\phi_a)G_2(\phi_a)d\mi&=\sum_{\el\ge0}\int_{S_\el\setminus S_{\el+1}}
G_1(\cm_\ct\phi_a)G_2(\phi_a)d\mi\ge\,\text{(due to (\ref{eq3.5}))} \nonumber \\
&\ge\sum_{\el\ge0}G_1(\thi_\el)\int_{S_\el\setminus S_{\el+1}}G_2(\phi_a)d\mi.  \label{eq3.6}
\end{align}
By the construction now of $\phi_a$ we note that
\[
\bigg(\phi_a/S_\el\setminus S_{\el+1}\bigg)^\ast=\bigg(g/((1-a)^{\el+1},(1-a)^\el]\bigg)^\ast,
\]
so (\ref{eq3.6}) becomes
\begin{align}
\int_XG(\cm_\ct\phi_a)G_2(\phi_a)d\mi \nonumber
\ge\sum_{\el\ge0}G_1\bigg(\frac{1}{(1-a)^\el}\int^{(1-a)^\el}_0g(u)du\bigg)\cdot\int^{(1-a)^\el}_{(1-a)^{\el+1}}G_2(g(u))du \nonumber\\
\ge\sum_{\el\ge0}G_1\bigg(\frac{1}{(1-a)^\el}\int^{(1-a)^\el}_0g(u)du\bigg)a(1-a)^\el G_2(g((1-a)^\el)) \nonumber \\
=\sum_{\el\ge0}G_1\bigg(\frac{1}{(1-a)^\el}\int^{(1-a)^\el}_0g(u)du\bigg)G_2
(g((1-a)^\el))|((1-a)^{\el+1},(1-a)^\el]|.  \label{eq3.7}
\end{align}
The sum in (\ref{eq3.7}) is a Riemman sum of the integral $\intl^1_0G_1\Big(\dfrac{1}{t}\intl_0^tg\Big)G_2(g(t))dt$, so as $a\ra0^+$,
we see that we have the needed inequality. The general case of the sharpness of Lemma \ref{lem3.1} for any $k$ can be proved along the same lines, integrating  $G_1(\cm_\ct\phi_a)\cdot G_2(\phi_a)$ on $S_{m_a}$ for each $a$, where $m_a\in\N$ is such that $(1-a)^{m_a+1}<k\le(1-a)^{m_a}$,and thus $(1-a)^{m_a}\ra k$, so by continuity reasons we have the result.

\section{The Bellman function}\label{sec4}
\noindent

We consider now a non-increasing function $g:(0,1]\ra\R^+$ and the quantities
\[
\begin{array}{l}
  v_g(L)=\intl^1_{t=0}\max\bigg(\dfrac{1}{t}\intl^t_0g,L\bigg)^pdt \ \ \text{and} \\ [2ex]
  u_g(L)=\intl^1_{t=0}g(t)\max\bigg(\dfrac{1}{t}\intl^t_0g,L\bigg)^{p-1}dt. \end{array}
\]
where $L\ge f$. We will prove the following
\setcounter{lem}{0}
\begin{lem}\label{lem4.1}
With the above notation the following equality holds for every $g:(0,1]\ra\R^+$,
\setcounter{equation}{0}
\begin{eqnarray}
v_g(L)=L^p-\frac{p}{p-1}fL^{p-1}+\frac{p}{p-1}u_g(L). \label{eq4.1}
\end{eqnarray}
\end{lem}
\begin{Proof}
We have that
\begin{align*}
v_g(L)&=\int^L_{\la=0}+\int^{+\infty}_{\la=L}p\la^{p-1}\bigg|\bigg\{t\in(0,1]:\max\bigg(
\frac{1}{t}\int^t_0g,L\bigg)\ge\la\bigg\}\bigg|d\la \\
&=L^p+\int^{+\infty}_{\la=L}p\la^{p-1}\bigg|\bigg\{t\in(0,1]:\frac{1}{t}\int^t_0
g\ge\la\bigg\}\bigg|d\la.
\end{align*}
We consider now for each $\la>L\ge f$, the unique $\bi(\la)\in(0,1]$ such that $\dfrac{1}{\bi(\la)}\intl^{\bi(\la)}_0g(u)du=\la$ (we suppose that $g(0^+)=+\infty$, without loss of the generality). Therefore,
\[
v_g(L)=L^p+\int^{+\infty}_{\la=L}p\la^{p-1}|A_\la|d\la,
\]
where
\[
A_\la=\bigg\{t\in(0,1]:\frac{1}{t}\int^t_0g>\la\bigg\}=(0,\bi(\la)). \ \ \text{So}
\]
\begin{align*}
v_g(L)&=L^p+\int^{+\infty}_{\la=L}p\la^{p-1}\bi(\la)d\la \\
&=L^p+\int^{+\infty}_{\la=L}p\la^{p-1}\bigg(\frac{1}{\la}\int^{\bi(\la)}_0g(u)du\bigg)d\la \\
&=L^p+\int^{+\infty}_{\la=L}p\la^{p-2}\bigg(\int_{\{u:\frac{1}{u}\intl^u_0g>\la\}}
g(u)du\bigg)d\la \\
&=L^p+\int^{+\infty}_{\la=L}p\la^{p-2}\bigg(\int_{\{u:\max\big(\frac{1}{u}\intl^u_0g,L\big)>\la\}}
g(u)du\bigg)d\la\\
&=L^p+\int^1_0g(t)\frac{p}{p-1}\big[\la^{p-1}\big]^{\max\big(\frac{1}{t}\intl^t_0g,L\big)}_{\la=L}dt\\
&=L^p-\frac{p}{p-1}L^{p-1}f+\frac{p}{p-1}u_g(L),
\end{align*}
where in the previous to the last inequality we have used Fubini's theorem.
Lemma 4.1 is now proved. \hs
\end{Proof}

We now prove the following
\begin{lem}\label{lem4.2}
For every $f$ and $F$ such that $0<f^p\le F$ and $L\ge f$ we have that
\[
\int_X\max(\cm_\ct\phi,L)^pd\mi\le\left\{\begin{array}{l}
                                           F\oo_p\Big(\dfrac{pL^{p-1}f-(p-1)L^p}{F}\Big)^p{F}, \;\text{if}\; L<\dfrac{p}{p-1}f \\  [2ex]
                                           F^p+\Big(\dfrac{p}{p-1}\Big)^p(F-f^p), \;\text{if}\;L\ge\dfrac{p}{p-1}f \end{array}\right.
\]
for every $\phi$ such that, $\intl_X\phi d\mi=f$ and $\intl_X\phi^pd\mi=F$.
\end{lem}
\begin{Proof}
We set $I=\intl_X\max(\cm_\ct\phi,L)^pd\mi$. Then
\begin{align*}
I&=\int^{+\infty}_{\la=0}p\la^{p-1}\mi(\{x\in X:\max(\cm_\ct\phi(x),L)>\la\})d\la \\
&=\int^L_{\la=0}+\int^{+\infty}_{\la=L}p\la^{p-1}\mi(\{x\in X:\max(\cm_\ct\phi(x),L)>\la\})d\la \\
&=II+III, \ \ \text{where}
\end{align*}
\[
II=\int^L_{\la=0}p\la^{p-1}d\la=L^p,
\]
since $(X,\mi)$ is a probability space, and
\[
III=\int^{+\infty}_{\la=L}p\la^{p-1}\mi(\{x\in X:\cm_\ct\phi(x)>\la\})d\la.
\]
By the weak type inequality (\ref{eq1.2}) we obtain that
\begin{align}
III&\le\int^{+\infty}_{\la=L}p\la^{p-1}\bigg(\frac{1}{\la}\int_{\{\cm_\ct\phi>\la\}}\phi d\mi\bigg)d\la \nonumber \\
&=\int^{+\infty}_{\la=L}p\la^{p-2}\bigg(\int_{\{\max(\cm_\ct\phi,L)>\la\}}\phi d\mi\bigg)d\la \nonumber \\
&=\int_X\phi(x)\bigg(\int^{\max(\cm_\ct\phi(x),L)}_{\la=L}p\la^{p-2}d\la\bigg)d\mi(x) \nonumber \\
&=\int_X\phi(x)\frac{p}{p-1}\big[\la^{p-1}\big]^{\max(\cm_\ct\phi(x),L)}_{\la=L}d\mi (x) \nonumber \\
&=\frac{p}{p-1}\int_X\phi(x)\max(\cm_\ct\phi(x),L)^{p-1}d\mi(x)-\frac{p}{p-1}L^{p-1}f. \label{eq4.2}
\end{align}
By (\ref{eq4.2}) then
\begin{align*}
III&\le\frac{p}{p-1}\bigg(\int_X\phi^pd\mi\bigg)^{1/p}\cdot\bigg(\int_X\max(\cm_\ct\phi,L)^p
\bigg)^{(p-1)/p}-\frac{p}{p-1}L^{p-1}f\Rightarrow \\
I&\le\frac{p}{p-1}F^{1/p}I^{(p-1)/p}+L^p-\frac{p}{p-1}L^{p-1}f\Rightarrow \\
\frac{I}{F}&\le\frac{p}{p-1}\bigg(\frac{I}{F}\bigg)^{(p-1)/p}+\frac{L^p-\dfrac{p}{p-1}L^{p-1}f}{F}
\Rightarrow
\end{align*}
\[
\Rightarrow pw^{p-1}-(p-1)w^p\ge\frac{pL^{p-1}f-(p-1)L^p}{F},
\]
where $w=\Big(\dfrac{I}{F}\Big)^{1/p}$. This gives
\begin{eqnarray}
-(p-1)w^p+pw^{p-1}=H_p(w)\ge\frac{pL^{p-1}f-(p-1)L^p}{F}  \label{eq4.3}
\end{eqnarray}
where the function $H_p$ is defined on $\Big[1,\dfrac{p}{p-1}\Big]$ with values on $[0,1]$.

We consider the function $h:[f,+\infty)\ra\R$ defined by
\[
h(t)=pt^{p-1}f-(p-1)t^p, \ \ t\ge f.
\]
Then
\begin{align*}
h'(t)&=p(p-1)t^{p-2}f-p(p-1)t^{p-1} \\
&=p(p-1)(f-t)t^{p-2}<0\Rightarrow\;\text{$h$ is strictly decreasing in it's domain}
\end{align*}
Therefore, $h(t)\le h(f)=f^p$ for every $t\ge f$, thus the right side of (\ref{eq4.3}) which we denote by $b$, is less than $f^p/F\le1$.

We consider two cases

i) $b\ge0$. Then we have that $b\in[0,1]$ and $H_p(\oo)\ge b$. If $w\le1$ then we must have that $I\le F$ which gives in view of the fact that $\oo_p(b)>1$, the inequality $I\le F[\oo_p(b)]^p$, that is our result. We consider now the case $w>1$.
Then since $H_p:\Big[1,\dfrac{p}{p-1}\Big]\ra[0,1]$ is strictly decreasing we have that
\begin{align*}
H_p(w)\ge b&\Rightarrow w\le\oo_p(b)\Rightarrow\frac{I}{F}\le[\oo_p(b)]^p \\
&\Rightarrow I\le F\oo_p\bigg(\frac{pL^{p-1}f-(p-1)L^p}{F}\bigg)^p,
\end{align*}
We have proved our Lemma in the first case.

ii) We consider now the second case: $b<0$ that is $L>L_0=\dfrac{p}{p-1}f$. Then
\[
I=\int_X\max(\cm_\ct\phi,L)^pd\mi=L^p+III
\]
where as we have seen
\begin{eqnarray}
III\le\int^{+\infty}_{\la=L}p\la^{p-2}\bigg(\int_{\{\cm_\ct\phi>\la\}}\phi d\mi\bigg)d\la. \label{eq4.4}
\end{eqnarray}
Since $L>L_0$ we conclude by (\ref{eq4.3}) that
\[
III\le\int^{+\infty}_{\la=L_0}p\la^{p-2}\bigg(\int_{\{\cm_\ct\phi>\la\}}\phi d\mi\bigg)d\la
=\int_X\max(\cm_\ct\phi,L_0)^pd\mi-L^p_0.
\]
By the case $L_0=\dfrac{p}{p-1}f$, which was treated in i) we conclude
\begin{align*}
\int_X\max(\cm_\ct\phi,L_0)^pd\mi&\le F\oo_p\bigg(\frac{pL^{p-1}_0f-(p-1)L^p_0}{F}\bigg)^p \\
&=F[\oo_p(0)]^p=F\bigg(\frac{p}{p-1}\bigg)^p.
\end{align*}
The above imply that
\[
I\le L^p+F\bigg(\frac{p}{p-1}\bigg)^p-L^p_0=L^p+\bigg(\frac{p}{p-1}\bigg)^p(F-f^p),
\]
which is our result in the second case.
Lemma \ref{lem4.2} is now proved. \hs
\end{Proof}

\section{Sharpness of Lemma \ref{lem4.1}}\label{sec5}
\noindent
We suppose now that $L<\dfrac{p}{p-1}f$ and look at the relations (\ref{eq4.1}) and (\ref{eq4.4}).
The first one is an inequality and states that
\setcounter{equation}{0}
\begin{eqnarray}
\int_X\max(\cm_\ct\phi,L)^pd\mi\le L^p-\frac{p}{p-1}L^{p-1}f+\frac{p}{p-1}
\int_X\phi\max(\cm_\ct\phi,L)^{p-1}d\mi\hspace*{-2cm}  \label{eq5.1}
\end{eqnarray}
while the second is an equality stating
\begin{eqnarray}
\int^1_0\max\bigg(\frac{1}{t}\int^t_0g,L\bigg)^pdt=L^p-\frac{p}{p-1}L^{p-1}f+\frac{p}{p-1}
\int^1_0g(t)\max\bigg(\frac{1}{t}\int^t_0g,L\bigg)^pdt.\hspace*{-2cm}  \label{eq5.2}
\end{eqnarray}
We fix $g:(0,1]\ra\R^+$. By Theorem 1 for
\[
\begin{array}{l}
  G_1(t)=\max(t,L)^p, \ \ t\ge0 \\ [2ex]
  G_2(t)=1, \ \ \text{and} \ \ k=1
\end{array}
\]
we have that
\[
\sup_{\phi^\ast=g}\int_X\max(\cm_\ct\phi,L)^pd\mi=v_g(L)
\]
while for
\[
\begin{array}{l}
  G_1(t)=\max(t,L)^{p-1}, \ \ t\ge0 \\ [2ex]
  G_2(t)=t, \ \ \text{and} \ \ k=1
\end{array}
\]
we see that
\[
\sup_{\phi^\ast=g}\int_X\phi\max(\cm_\ct\phi,L)^{p-1}d\mi=u_g(L).
\]
That is if we leave the $\phi$'s to move along the rearrangements of $g$ in (\ref{eq4.1}) we produce the equality (\ref{eq4.4}). During the
proof of Lemma 4.1 we have also used the following inequality
\begin{eqnarray}
\int_X\phi\max(\cm_\ct\phi,L)^{p-1}d\mi\le\bigg(\int_X\phi^pd\mi\bigg)^{1/p}
\bigg(\int_X\max(\cm_\ct\phi,L)^pd\mi\bigg)^{p-1/p}.\hspace*{-1cm}  \label{eq5.3}
\end{eqnarray}
For the proof of Lemma 4.2 we used inequalities only in the above two mentioned points. The first is attained if we use (\ref{eq4.1}) and the discussion before. For the second we conclude that we need to find a sequence $g_n:(0,1]\ra\R^+$ with $\intl^1_0g_n(u)du=f$ and $\intl^1_0g^p_n(u)du=F$ for which
\[
\int^1_0g_n(t)\max\bigg(\frac{1}{t}\int^t_0g_n,L\bigg)^{p-1}dt\approx\bigg(
\int^1_0g^p_n\bigg)^{1/p}\cdot\bigg(\int^1_0\max\bigg(\frac{1}{t}\int^t_0g_n,L\bigg)dt\bigg)^{(p-1)/p}
\]
that is we need equality in a Holder inequality.
Therefore, we are forced to search for a $g:(0,1]\ra\R^+$ with
\[
\int^1_0g(u)du=f \ \ \text{and} \ \ \int^1_0g^p(u)du=F
\]
for which
\begin{eqnarray}
\max\bigg(\frac{1}{t}\int^t_0g,L\bigg)=cg(t),  \ \ \text{for} \ \ t\in(0,1] \label{eq5.4}
\end{eqnarray}
where
\[
c=\oo_p\bigg(\frac{pL^{p-1}f-(p-1)L^p}{F}\bigg).
\]
We state it as
\setcounter{lem}{0}
\begin{lem}\label{lem5.1}
There exists $g:(0,1]\ra\R^+$ non-increasing, continuous for which the above three equations for the constants $f,F$ and $c$ hold, in case where $L<\dfrac{p}{p-1}f$.
\end{lem}
\begin{Proof}
We set
\begin{eqnarray}
g(t)=\left\{\begin{array}{ccc}
                    Kt^{-1+\frac{1}{c}}, & \text{if} & t\in[0,\ga] \\ [2.5ex]
                    \frac{L}{c}, & \text{if} & t\in[\ga,1]
                  \end{array}\right.
\end{eqnarray}
where $\ga$ and $K$ are such that $\dfrac{1}{\ga}\intl^\ga_0g(u)du=L$, that is
\begin{eqnarray}
Kc\ga^{-1+\frac{1}{c}}=L.  \label{eq5.5}
\end{eqnarray}
It is obvious that $g$ is continuous, non-increasing and satisfies (\ref{eq5.4}). We are going to find now the constant $\ga$ in a way that
\[
\int^1_0g^p(u)du=F\Leftrightarrow\frac{K^p\big[t^{-p+\frac{p}{c}+1}\big]^\ga_{t=0}}
{\Big(-p+\dfrac{p}{c}+1\Big)}+\frac{L^p}{c^p}(1-\ga)=F\Leftrightarrow
\]
\begin{align}
&\frac{K^pc^p\ga^{-p+\frac{p}{c}+1}}{c^p\Big(-p+\dfrac{p}{c}+1\Big)}+\frac{L^p}{c^p}
(1-\ga)=F \nonumber \\
&\Leftrightarrow\frac{c^pK^p\ga^{-p+\frac{p}{c}+1}}{-(p-1)c^p+pc^{p-1}}+
\frac{L^p}{c^p}(1-\ga)=F  \label{eq5.6}.
\end{align}
Since (\ref{eq5.5}) holds (\ref{eq5.6}) becomes
\begin{eqnarray}
\frac{L^p\cdot\ga}{-(p-1)c^p+pc^{p-1}}+\frac{L^p}{c^p}(1-\ga)=F.  \label{eq5.7}
\end{eqnarray}
By the definition of $c$ we have that
\[
-(p-1)c^p+pc^{p-1}=\frac{pL^{p-1}f-(p-1)L^p}{F}=b,
\]
so (\ref{eq5.7}) becomes
\[
\frac{FL^p\cdot\ga}{pL^{p-1}f-(p-1)L^p}+\frac{L^p}{c^p}(1-\ga)=F\;\Leftrightarrow
\]
\[
\Leftrightarrow\;\ga=\frac{F-L^p/c^p}{L^p\Big(\dfrac{1}{b}-\dfrac{1}{c^p}\Big)}.
\]
We need to see that $\ga\in[0,1]$. Obviously we have that
\[
L^p\le\int_X\max(\cm_\ct\phi,L)^pd\mi
\]
for any $\phi$ such that $\intl_X\phi d\mi=f$ and $\intl_X\phi^pd\mi=F$. Additionally
\[
\int_X\max(\cm_\ct\phi,L)^pd\mi\le[\oo_p(b)]^p\cdot F=c^pF\Rightarrow F-L^p/c^p\ge0.
\]
Further $c$ satisfies $-(p-1)c^p+pc^{p-1}=b$ as it is mentioned before thus $p(c^p-c^{p-1})=c^p-b\Rightarrow c^p-b>0\Rightarrow\dfrac{1}{b}-\dfrac{1}{c^p}>0$.
From the above two inequalities we see that $\ga\ge0$
We prove now that $\ga\le1\Leftrightarrow$
\[
F-\frac{L^p}{c^p}\le\frac{L^p}{b}-\frac{L^p}{c^p}\Leftrightarrow
\]
\[
F\cdot b\le L^p\Leftrightarrow F\cdot\frac{pL^{p-1}f-(p-1)L^p}{F}\le L^p\Leftrightarrow L^{p-1}f\le L^p,
\]
which is true because of the fact that always $L\ge f$.

We consider now the function $g$ as defined before with
\[
\ga=\frac{F-L^p/c^p}{L^p\Big(\dfrac{1}{b}-\dfrac{1}{c^p}\Big)}\in[0,1].
\]
We prove that we additionally have that
\begin{align*}
\int^1_0g(u)du=f&\Leftrightarrow\int^\ga_0Kt^{-1+\frac{1}{c}}dt+\frac{L}{c}(1-\ga)=f \\
&\Leftrightarrow Kc\ga^{1/c}+\frac{L}{c}(1-\ga)=f \\
&\Leftrightarrow\;\text{(since $Kc=L\ga^{1-\frac{1}{c}}$)}
\end{align*}
\[
L\ga+\frac{L}{c}(1-\ga)=f\Leftrightarrow\ga=\frac{f-L/c}{L\Big(1-\dfrac{1}{c}\Big)},
\]
So we need to check that
\begin{align*}
\frac{f-\dfrac{L}{c}}{L\Big(1-\dfrac{1}{c}\Big)}&=\frac{F-\dfrac{L^p}{c^p}}
{L^p\Big(\dfrac{1}{b}-\dfrac{1}{c^p}\Big)}\Leftrightarrow \\
\frac{fc-L}{(c-1)}&=\frac{Fc^p-L^p}{L^{p-1}\Big(\dfrac{c^p}{b}-1\Big)}\Leftrightarrow
\end{align*}
\begin{eqnarray}
b=\frac{c^{p-1}(fc-L)L^{p-1}}{F(c^p-c^{p-1})-L^p+fL^{p-1}}.  \label{eq5.8}
\end{eqnarray}
Because now of the relation
\[
c^p-c^{p-1}=\frac{-b+c^p}{p},
\]
(\ref{eq5.8}) becomes
\begin{eqnarray}
b=\frac{c^{p-1}(fc-L)L^{p-1}}{\dfrac{F}{p}(-b+c^p)-L^p+fL^{p-1}}.  \label{eq5.9}
\end{eqnarray}
On the other hand
\begin{align*}
\frac{F}{p}(-b+c^p)-L^p+fL^{p-1}&=\frac{F}{p}\bigg(-\frac{pL^{p-1}f-(p-1)L^p}{F}+c^p\bigg)-L^p+fL^{p-1} \\
&=-L^{p-1}f+\frac{p-1}{p}L^p+\frac{F}{p}c^p-L^p+fL^{p-1} \\
&=\frac{F}{p}c^p-\frac{L^p}{p}=\frac{Fc^p-L^p}{p}.
\end{align*}
Thus (\ref{eq5.9}) is equivalent to
\[
b=\frac{pc^{p-1}(fc-L)L^{p-1}}{Fc^p-L^p}\;\Leftrightarrow \]
\begin{align*}
&\Leftrightarrow\frac{pc^pf}{L}-pc^{p-1}=b\bigg(\frac{Fc^p}{L^p}-1\bigg)\Leftrightarrow \text{(since $pc^{p-1}=b+(p-1)c^p)$}\\
&\Leftrightarrow\frac{pc^pf}{L}-b-(p-1)c^p=bF\frac{c^p}{L^p}-b\Leftrightarrow
\frac{pf}{L}-(p-1)=b\frac{F}{L^p}\Leftrightarrow
\end{align*}
\[
b=\frac{pL^{p-1}f-(p-1)L^p}{F}
\]
which is true from the definition of $b$.

That is we derived Lemma \ref{lem5.1}.  \hs
\end{Proof}
We turn now to the case $L\ge\dfrac{p}{p-1}f$.
For this one we need to construct a sequence $(g_n)_n$ with $g_n:(0,1]\ra\R^+$ non-increasing and continuous such that
\[
\int^1_0g_n(u)du=f, \ \ \int^1_0g^p_n(u)du=F \ \ \text{and}
\]
\[
\lim_n\int^1_0\max\bigg(\frac{1}{t}\int^t_0g_n,L\bigg)^pdt\ge L^p+\bigg(\frac{p}{p-1}\bigg)^p(F-f^p)
\]
where $L\ge\dfrac{p}{p-1}f$.

We set as before
\[
g_n(t)=\left\{\begin{array}{ll}
              k_nt^{-1+\frac{1}{c_n}}, & t\in(0,\ga_n] \\ [1.5ex]
              \dfrac{L_n}{c}, & t\in[\ga_n,1]
            \end{array}\right.
\]
where $L_n\nearrow L_0=\dfrac{p}{p-1}f$,
\[
\ga_n=\frac{F-L^p_n/c^p_n}{L^p_n\Big(\dfrac{1}{b_n}-\dfrac{1}{c^p_n}\Big)}=\frac{f-L_n/c_n}{L_n
\Big(1-\dfrac{1}{c_n}\Big)}
\]
where $c_n=\oo_p(b_n)$, $b_n=\dfrac{pL^{p-1}_nf-(p-1)L^p_n}{F}$ and $k_n$ is such that $k_nc_n\ga_n^{-1+\frac{1}{c_n}}=L_n$.
Since $L_n\ra L_0$ we have that $b_n\ra0$, $c_n\ra\dfrac{p}{p-1}$ and $\ga_n\searrow\dfrac{f-L_0\dfrac{p-1}{p}}{L_0\Big(1-\dfrac{p}{p-1}\Big)}=0$.
According to the first case (where $L<\dfrac{p}{p-1}f$) we have that
\[
\int^1_0\max\bigg(\frac{1}{t}\int^t_0g_n,L_n\bigg)^pdt=[\oo_p(b_n)]^pF\ra\bigg(
\frac{p}{p-1}\bigg)^pF.
\]
Now for $L\ge\dfrac{p}{p-1}f$,
\begin{align}
\int^1_0\max\bigg(\frac{1}{t}\int^t_0g_n,L\bigg)^pdt=&\,L^p+\int^{+\infty}_{\la=L}
p\la^{p-2}\bigg(\int_{\big\{u:\frac{1}{u}\int^u_0g_n>\la\big\}}g_n(u)du\bigg)d\la\nonumber \\
\overset{L>L_0}{=}&L^p+\int^{+\infty}_{\la=L_0}p\la^{p-2}\bigg(
\int_{\big\{u:\frac{1}{u}\intl^t_0g_n>\la\big\}}g_n(u)du\bigg)d\la\nonumber\\
&-\int^L_{\la=L_0}p\la^{p-2}\bigg(\int_{\big\{u:\frac{1}{u}\intl^u_0g_n>\la\big\}}
g_n(u)du\bigg)d\la \nonumber \\
=&L^p-L^p_0+\int^1_0\max\bigg(\frac{1}{t}\int^t_0g_n,L_0\bigg)^pdt \nonumber \\
&-\int^L_{\la=L_0}p\la^{p-2}\bigg(\int_{\big\{u:\frac{1}{u}\intl^u_0g_n>\la\big\}}
g_n(u)du\bigg)d\la  \label{eq5.10}
\end{align}
By definition of the functions $g_n$ we have that
\[
\max\bigg(\frac{1}{t}\int^t_0g_n,L_n\bigg)=\oo_p(b_n)g_n(t).
\]
Thus
\begin{align*}
\int^1_0\max\bigg(\frac{1}{t}\int^t_0g_n,L_0\bigg)^pdt&\ge\int^1_0\max
\bigg(\frac{1}{t}\int^t_0g_n,L_n\bigg)^pdt \\
&=[\oo_p(b_n)]^p\int^1_0g^p_n(u)du=F[\oo_p(b_n)]^p, \ \ \text{for every $n$}
\end{align*}
and so
\[
\lim_n\int^1_0\max\bigg(\frac{1}{t}\int^t_0g_n,L_0\bigg)^pdt=F\bigg(\frac{p}{p-1}\bigg)^p.
\]
At last
\[
a_n(L)=\int^L_{\la=L_0}p\la^{p-2}\bigg(\int_{\big\{t:\frac{1}{t}\intl^t_0g_n>\la\big\}}
g_n(u)du\bigg)d\la
\]
satisfies for a given $L\ge L_0$
\begin{align}
a_n(L)&\le\int^L_{\la=L_0}p\la^{p-2}\bigg(\int_{\big\{t:\frac{1}{t}\intl^t_0g_n>L_0\big\}}
g_n(u)du\bigg)d\la \nonumber \\
&=\bigg(\int_{\big\{t:\frac{1}{t}\intl^t_0g_n>L_0\big\}}g_n(u)du\bigg)
\int^L_{\la=L_0}p\la^{p-2}d\la\nonumber \\
&=\ti_L\cdot\int_{\big\{t:\frac{1}{t}\intl^t_0g_n>L_0\big\}}g_n(u)du.  \label{eq5.11}
\end{align}
Note then that
\[
\bigg|\bigg\{t\in(0,1]:\frac{1}{t}\int^t_0g_n\ge L_0\bigg\}\bigg|\le
\bigg|\bigg\{t\in(0,1]:\frac{1}{t}\int^t_0g_n\ge L_n\bigg\}\bigg|=\ga_n,
\]
because $\ga_n$ is the unique element of $(0,1]$ such that $\dfrac{1}{\ga_n}\intl^{\ga_n}_0g_n=L_n$.

Since $\ga_n\ra0$, from (\ref{eq5.11}) we deduce that $a_n(L)\ra0$, as $n\ra\infty$, thus from (\ref{eq5.10})
\[
\lim_n\int^1_0\max\bigg(\frac{1}{t}\int^t_0g_n,L\bigg)^pdt\ge L^p-L^p_0
+\bigg(\frac{p}{p-1}\bigg)^pF=L^p+\bigg(\frac{p}{p-1}\bigg)^p(F-f^p),
\]
which is the result we needed to prove.
From Lemma \ref{lem5.1} and the calculations after it's proof we conclude the sharpness of Lemma \ref{lem4.1}.

\section{Conclusions}\label{sec6}
By providing a generalization of the symmetrization principle given in \cite{10} we give another
proof of the computation for the Bellman function of three variables of the dyadic maximal operator,different from those given in \cite{3} and \cite{11}.

\vspace*{2cm}
\noindent
{\bf E.Nikolidakis, A. Melas:} Department of Mathematics \vspace*{0.1cm} \\
National and Kapodistrian University of Athens, \vspace*{0.1cm} \\
GR-157 84 Panepistimioupolis, \vspace*{0.1cm} \\
Athens, Greece \vspace*{0.1cm} \\
E-mail addresses: lefteris@math.uoc.gr, amelas@math.uoa.gr \\

\end{document}